\documentclass[12pt,reqno]{amsart}
\usepackage{amsmath}
\usepackage{amsthm}
\usepackage{mathrsfs}
\usepackage{amsfonts} 
\usepackage{epsfig}

\numberwithin{equation}{section}

\allowdisplaybreaks

\newtheorem{theorem}{Theorem}[section]

\newtheorem{proposition}[theorem]{Proposition}
\newtheorem{lemma}[theorem]{Lemma}

\newtheorem{definition}[theorem]{Definition}

\renewcommand{\epsilon}{\varepsilon}

\newcommand{\R}{\mathbb{R}}
\newcommand{\N}{\mathbb{N}}
\newcommand{\Z}{\mathbb{Z}}
 
\newcommand{\PR}{\mathbb{P}}

\newcommand{\good}{\mathcal{U}}
\newcommand{\parent}[1]{{#1}^+}
\newcommand{\child}[1]{{#1}^-}
\newcommand{\zreg}{{\rm RG}}
\newcommand{\zstrreg}{{\rm SRG}(Z)}
\newcommand{\ladder}{\mathcal{B}}
\newcommand{\degr}[1]{{\rm deg}(#1)}

\newcommand{\kappan}{T/2}
\newcommand{\tminusk}{T/2}
\newcommand{\degg}{d}

\newcommand{\cbeta}{c_{1,T}}
\newcommand{\ctbeta}{c_{2,T}}
\newcommand{\cbetastar}{c^*_{1,T}}
\newcommand{\ctbetastar}{c^*_{2,T}}
\newcommand{\tdegg}{\mathcal{T}_\degg}
\newcommand{\tcritone}{T_c^1}
\newcommand{\tcrittwo}{T_c^2}
\newcommand{\hnew}{\eta}
\newcommand{\tint}{\mathscr{T}}
\newcommand{\coupling}{\mathcal{C}}
\newcommand{\tree}{\mathcal{T}}
\newcommand{\treeg}{G}
\newcommand{\graph}{G}

\newcommand{\ubl}{{\rm UnTouch}}

\newcommand{\de}{d_0}
\newcommand{\dei}{d_0^{-1}}
\newcommand{\deit}{d_0^{-2}}
\newcommand{\deith}{d_0^{-3}}

\newcommand{\y}{Y}

\newcommand{\forced}{{\rm Found}}

\newcommand{\laddergiven}{\mathcal{B}_0}

\newcommand{\dd}{\, \mathrm{d}}

\title[Infinite cycles in the random stirring model]{Infinite cycles in \\ the random stirring model on trees
}
\date{}

\author[A.~Hammond]{Alan HAMMOND}
\address{ Department of Statistics, 
           University of Oxford,
                  1 South Parks Road,
                  Oxford, OX1 3TG, U.K. } 
\email{hammond@stats.ox.ac.uk}

\thanks{Department of Statistics, Oxford University. Supported by  U.K. EPSRC grant EP/I004378/1 and U.S. NSF grant OISE-07-30136.} 

\begin{document}

\begin{abstract}
We prove that, in the random stirring model 
of parameter $T > 0$ on an infinite rooted tree each of whose vertices has at least two offspring, infinite cycles exist almost surely, provided that $T$ is sufficiently high.
\end{abstract}

\maketitle

\section{Introduction}

It has long been recognized in the physics literature that ``spatial random permutations'' -- laws on permutations defined from spatial models -- are intimately related to the behaviour of low-temperature gases. We begin by recounting briefly the first such example.

In 1953, Feynman \cite{feynman} wrote the quantum-mechanical partition function for helium as a sum over the energy associated to certain interacting Brownian particles that may interchange their positions over a finite-time interval. He argued that the $\lambda$-transition undergone by the gas at low temperature is reflected by the appearance of large cycles in a measure on permutations naturally associated to this representation of the partition function. 

 To describe more precisely this law on permutations, we consider only the hard-core instance  (which formally corresponds to a potential that is $+\infty$ when it is non-zero; the potential for helium is finite, varying from highly positive below the atomic radius to slightly negative on the order of this radius).
Fix dimension $d \geq 2$, as well as a ``time" (or inverse-temperature) parameter $\beta \in (0,\infty)$ and a small ``interaction-range" parameter $r > 0$.
Scatter a large number $N$ of points independently and uniformly in the $d$-dimensional torus of volume $N$, and run from each of them an independent Brownian motion for time $\beta$; the law appearing in Feynman's repesentation of the gas is obtained by conditioning the system on the time-$\beta$ configuration of $N$ points coinciding with the time-zero configuration, and on the avoidance constraint that  no pair of points be at distance less than $r$ at any time $t \in [0,\beta]$. A random permutation arises by mapping each point at time-zero to the point at time-$\beta$ obtained by following the Brownian path beginning at the point during the period $[0,\beta]$. 

It is anticipated that, when $d \geq 3$, and $r > 0$ is fixed at a small enough value,  cycles of macroscopic volume (of order $N$) appear in the model, provided that $\beta$ exceeds a critical value, and that this behaviour reflects the condensation of the gas at low-temperature. 
 The behaviour of Feynman's model is understood rigorously only in the non-interacting case (where formally $r=0$, so that no avoidance conditioning is applied). Here, 
the existence of the critical value for the presence of macroscopic cycles
was proved by S{\"u}t\H{o} \cite{sutoone,sutotwo}, who showed that it coincides with the critical density of the ideal Bose gas identified by Einstein \cite{einstein}. Some extensions of these results to other non-interacting models are made in \cite{buspatial}. To the best of our knowledge, no direct argument has been made to establish the existence of large cycles in an interacting model in a Euclidean setting.  

It is a physically important and mathematically very interesting question, then,  to prove the presence of large cycles in natural models of spatial random permutations.

\subsection{Main results}
In this article, we study the cycles in the random stirring model on a tree.  The random stirring model on a given graph $G$
is the stochastic process $\sigma$ mapping $[0,\infty)$ to permutations of the vertex set of $G$ which starts at the identity and under which the transposition associated to each edge in $G$
is performed at each of the points in a Poisson process with mean one, independently for each edge. 
For each $T \in [0,\infty)$, we will refer to the marginal law $\sigma_T$ as the random stirring model with parameter~$T$.

Omer Angel \cite{angel} has proved that, on a regular tree of degree at least five, and for a certain interval of values of $T$, the random stirring model 
with parameter $T$ on the tree has infinite cycles almost surely.

We now state our main theorem.
 Our result
develops Angel's, by dispensing with the hypothesis that vertices have at least four offspring, and by being applicable for all sufficiently high $T$ (though his result begins to apply for slightly smaller values of $T$, as we shortly discuss). 
Terms from graph theory are reviewed in Subsection \ref{secctrw}.
\begin{theorem}\label{thmone}
Let $G$ be any infinite rooted  tree of uniformly bounded degree each of whose vertices has at least two offspring. Then there exists $T_0 \in (0,\infty)$ such that
if $T \geq T_0$ then the random stirring model with parameter $T$ contains infinite cycles almost surely.
\end{theorem}
The second theorem quantifies the value of $T_0$ for high-degree trees. 
\begin{theorem}\label{thmtwo}
Let $\degg \geq 55$. 
Let $G$ be an infinite rooted  tree of uniformly bounded degree 
each of whose  vertices has at least $\degg$ offspring.
Then we may choose $T_0 = 101 \degg^{-1}$ in the statement of Theorem~\ref{thmone}.
\end{theorem}

\subsection{Literature on the random stirring model}

The random stirring model (which is also called the random interchange model) was introduced in \cite{harris}. Its physical relevance  was indicated by B{\'a}lint T{\'o}th \cite{toth}, who used it to give a representation of the spin-$1/2$ Heisenberg ferromagnet; the lecture notes \cite{guw} contain an overview of this topic. Recent mathematical progress on the model 
includes the resolution of Aldous' conjecture identifying its spectral gap \cite{clr}, and a formula for the probability that the random permutation consists of a single cycle \cite{alonkozma}.

The emergence of a giant component under percolation on the complete graph as the percolation parameter 
$p$ increases through values near $1/n$ has been intensively studied. This transition is accompanied by the appearance of large-scale cycles in the associated random stirring model, as we now review.
Under the uniform measure on permutations on a finite set $V$, the lengths of cycles, normalized by $\vert V \vert$, and listed in decreasing order, converge to the Poisson-Dirichlet distribution with parameter one. Studying a model very closely related to the random stirring model for the complete graph, Oded Schramm considered the law on permutations of an $n$-point set given by composing $tn$ uniform random transpositions \cite{comprantrans}. (We will refer to this law as the $(n,t)$-random composition model.) Under this law, say that two vertices are connected if a transposition has been made on the edge between them.
Reflecting the emergence of a giant component in percolation on the $n$-point complete graph at parameter $p = 1/n$, the $(n,t)$-random composition model with $t = 1 + \epsilon$ contains a giant connected component of some density $\theta(\epsilon) \in (0,1)$. 
It is shown in \cite{comprantrans} that the ordered list of cycle lengths normalized by $\theta(\epsilon) n$ converges to the Poisson-Dirichlet distribution of parameter one; that is, a local equilibrium for large cycles inside the giant connected component is achieved as soon as this component becomes macroscopic.
Nathana{\"e}l Berestycki \cite{berestycki} has given a short proof that a cycle exists of size $\Theta(n)$ when $t = 1 + \epsilon$.
 
\subsection{The cyclic-time random walk}\label{secctrw}

Our analysis of the random stirring model exploits a closely related dependent random walk which B{\'a}lint T{\'o}th discusses in the proof of Theorem 1 in \cite{toth} and which Omer Angel in~\cite{angel}  calls the cyclic-time random walk. We now introduce further notation and define this walk.

We begin by recalling some graph-theoretic notation. The vertex and edge-sets of a given graph $G$ will be denoted by $V(G)$ and $E(G)$.
A graph is rooted if it has a distinguished vertex, the root, that we will denote by $\phi$. We write $d:V(G) \times V(G) \to \N$ for the graphical distance on $G$. A connected graph without cycles is called a tree. In a tree, there is a unique simple path $P_v$ leading from any given vertex $v$ to the root; the first element after $v$ on $P_v$ is called the parent of $v$, and each vertex is called an offspring of its parent. For $v,w \in V(G)$, $v$ is called a descendent of $w$ if $w$ is a vertex in $P_v$; $v$ is called a strict descendent of $w$ if, in addition, it is not equal to~$w$. Note that the set of descendents of a given vertex induces a subtree of $G$ (which we call the descendent tree of the vertex). 
  For each vertex $v \in V(G)$, we write $E_v$ for the set of edges incident to $v$; we write  $\degr{v} = \big\vert E_v \big\vert$ for the degree of~$v$.
 For each edge $e \in E(\graph)$, the incident vertex of $e$ closer to $\phi$ will be called the parent vertex of $e$ and will be denoted by $e^+$; the other, called the child vertex of $e$ and labelled $e^-$. 
  
Throughout we take $\graph$ to be a rooted tree whose vertex degree is uniformly bounded and each of whose vertices has at least two offspring. Sometimes we will further invoke the hypothesis that, for some given $\degg \geq 2$,
\begin{equation}\label{eqatleastdegg}
 \text{each vertex in $G$ has at least $\degg$ offspring.} 
\end{equation}

We now present a construction of cyclic-time random walk. Throughout, fix $T \in (0,\infty)$. For convenience, suppose that $\graph$ is embedded in $\R^2$, so that each element of $V(\graph)$ is identified with a point in $\R^2$ and each element $e \in E(\graph)$ with the line segment $[v_1,v_2] \subseteq \R^2$ where $e = (v_1,v_2)$ for $v_1,v_2 \in V(\graph)$. For each $v \in V(\graph)$, let the pole at $v$,
 $\{ v \} \times [0,T) \subseteq \R^3$, denote the line segment of length $T$ that rises vertically from $v$. Elements of $E(\graph) \times [0,T)$ will be called bars.
     The bar $b = (e,h)$ is said to be supported on the edge $e$ and to have height $h$.
   Note that the bar $(e,h)$ is a horizontal line segment which intersects the poles at $e^+$ and $e^-$; the intersection points $(e^+,h)$ and $(e^-,h)$ will be called the joints of $(e,h)$. 

The bar set $E(\graph) \times [0,T)$ carries the product of counting and Lebesgue measure on its components. (As a shorthand, we will refer to this product measure simply as Lebesgue measure.)

Let $(v,h) \in V(\graph) \times [0,T)$. Unit-speed cyclic upward motion from $(v,h)$ is the process $[0,\infty) \to \{ v \} \times [0,T): t \to \big(v, (h + t) \, \textrm{mod} \, T \big)$.

Let $\laddergiven \subseteq E(\graph) \times [0,T)$ be a collection of bars. Cyclic-time meander $X^{\laddergiven}_{(v,h)}:[0,\infty) \to V(\graph) \times [0,T)$,
 among $\laddergiven$ and with initial condition $(v,h) \in V(\graph) \times [0,T)$, is the following process. First, $X_{(v,h)}(0) = (v,h)$; the process pursues unit-speed cyclic upward motion from $(v,h)$  until (the possibly infinite time at which) it reaches the joint of a bar in $\laddergiven$, when it jumps to the other joint of this bar. The process $X^{\laddergiven}_{(v,h)}$ then continues by iterating the same rule, until it is defined on all of $[0,\infty)$. The process is chosen to be right-continuous with left limits. 
 We write $X^{\laddergiven}$ 
for  $X^{\laddergiven}_{(\phi,0)}$. (There are choices of $\laddergiven$ for which these rules fail to define $X^{\laddergiven}_{(v,h)}$ on all of $[0,\infty)$. It is a simple matter to verify that this difficulty does not arise in the case that is relevant to us and which we now discuss.)

We write $\PR_T$ for a probability measure carrying a bar collection $\ladder \subseteq E(\graph) \times [0,T)$ 
having Poisson law with intensity one with respect to Lebesgue measure. 
Cyclic-time random meander with parameter~$T$ is the random process $X^{\ladder}$.  We write $X$ in place of $X^\ladder$ and call $X$ in shorthand a meander.

\begin{figure}
\centering\epsfig{file=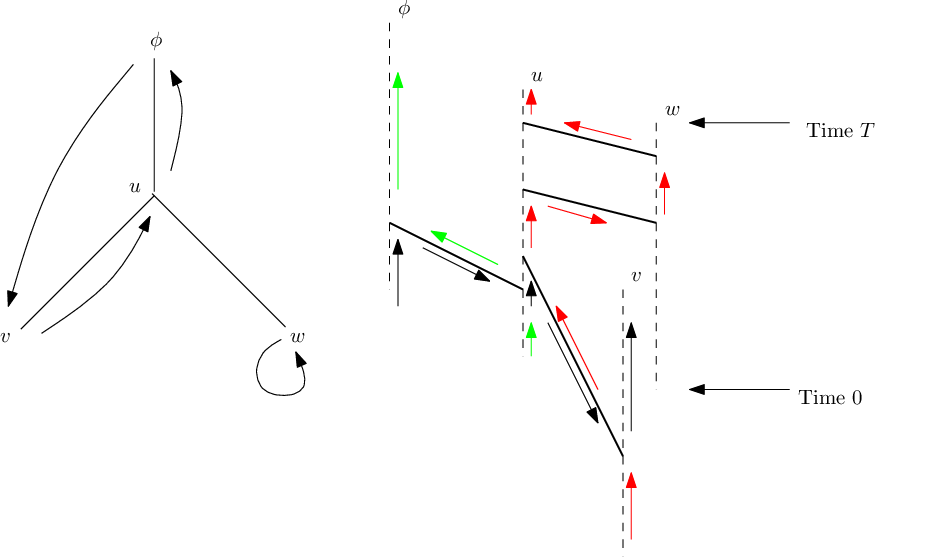, width=14cm}
\caption{For the graph shown on the left, cyclic-time random meander $X$ departing from $(\phi,0)$ is illustrated on the right. The right-hand sketch depicts a construction in $\R^3$ in which the poles associated to vertices are the vertical dashed lines and the bars in $\ladder$ are the horizontal black lines.
 Assume that there  are no bars in $\ladder$ supported on edges that  connect vertices $v$ and $w$ of $\phi$ to their offspring. The trajectory of the meander from $(\phi,0)$ is divided into three intervals of duration $T$, at the end of which, the meander returns to $(\phi,0)$. These three sub-trajectories are indicated in black, red and green in the right-hand sketch. As the left-hand sketch shows, the cycle of $\phi$ in the associated permutation 
thus has three elements.}\label{cycleexample}
\end{figure}

Cyclic-time random walk (begun at $v \in V(\graph)$ and with parameter~$T$) is the vertex-valued process given by projecting $X_{(v,0)}:[0,\infty) \to V(\graph) \times [0,T)$ onto $V(\graph)$. We denote it by $\y_v$ and write $Y$ in place of $Y_\phi$. 
    Note that the random stirring model with parameter~$T$ is the law of the random map $V(\graph) \to V(\graph): v \to Y_v(T)$.
See Figure~\ref{cycleexample} for an illustration.

Following Angel, we say that cyclic-time random walk $Y$ is transient
if there is positive probability that it never returns to the root, in the sense that 
 there exists $s_0 > 0$ such that $\phi \not\in Y(s_0,\infty)$. 
 
Theorem~\ref{thmone} will follow directly from the next proposition.
\begin{proposition}\label{propone}
Let $G$ be any infinite rooted  tree of uniformly bounded degree each of whose vertices has at least two offspring.
 Then there exists $T_0 \in (0,\infty)$ such that
if $T \geq T_0$ then cyclic-time random walk $Y$
is transient. 
\end{proposition}


Proposition~\ref{propone} and Theorem~\ref{thmone} are proved in Section \ref{secproofs}. Theorem~\ref{thmtwo} is obtained by reprising these arguments and developing quantitative counterparts to limiting assertions made along the way. Its proof is given in Appendix $A$.

\subsection{The sharp transition conjecture and trees of high degree}

For any given graph $G$ on which the random stirring model is well-defined, 
let $\tint^G$ denote the set of  $T > 0$ such that the random stirring model on $G$ with parameter $T$ has infinite cycles almost surely. 
Note that $T \not\in \tint^G$ unless the bond percolation on $G$ given by the set of edges that support a bar in $\ladder$ has an infinite component. As noted in \cite{angel}, this implies that 
$\big[ 0 , - \log ( 1 - p_c ) \big) \cap \tint^G = \emptyset$,
where $p_c = p_c(G)$ denotes the critical value for bond percolation on $G$. Writing $\tree_\degg$ for the rooted regular tree each of whose vertices has $\degg$ offspring, note that 
$p_c(\tree_\degg) = \degg^{-1}$, and thus that,
if $\degg \geq 2$, then 
\begin{equation}\label{eqdegglb} 
\big[ 0 , \degg^{-1} + \tfrac{1}{2}\degg^{-2} \big) \cap \tint^{\tree_\degg} = \emptyset \, .
\end{equation}
Define the critical points $\tcritone(G) = \inf \tint^G$ and $\tcrittwo(G) = \sup \, \big( [0,\infty) \setminus \tint^G \big)$. 
Note that
$\tcritone(G) \leq \tcrittwo(G)$ trivially. Conjecture $9$ of \cite{angel} claims that, for any graph $G$ for which $\tint^G$ is non-empty, these two critical points are equal. 
The present work and \cite{angel} go some way to verifying the conjecture for high-degree trees:  
we will prove the next result in Appendix $B$.
\begin{theorem}\label{thmthree} 
For any $\epsilon > 0$, there exists $d_0 \in \N$ such that if $\degg \geq d_0$
then  $\big[ \degg^{-1} + (\tfrac{7}{6} + \epsilon) \degg^{-2}, \infty) \subseteq \tint^{\tree_{\degg}}$. 
\end{theorem}
This deduction and (\ref{eqdegglb}) show that the discrepancy $\tcrittwo(G) - \tcritone(G)$ is $O(d^{-2})$ for high $\degg$.
In \cite{hammondtwo}, the discrepancy is shown to be zero for regular trees of high-degree, thus confirming Conjecture $9$ of \cite{angel} for such trees.


\vspace{2mm}

\noindent{\bf Acknowledgments.} I thank Christophe Garban and Daniel Ueltschi for useful discussions, and L\H{o}rinc S\'ark\'any and a referee for comments on the manuscript.

\section{Proofs}\label{secproofs}

We now begin to define and describe the elements needed to prove Proposition~\ref{propone}, after which, we will give the  proof of this result; Theorem~\ref{thmone} will then be an immediate consequence.

\subsection{Preliminaries}

Here we record a simple observation regarding future of cyclic-time random meander given its past.
\begin{lemma}\label{lemubl}
Let $t > 0$.
 Consider the law $\PR_T$ given $X:[0,t] \to V(\treeg) \times [0,T)$. Let $\forced_t \subseteq E(\treeg) \times [0,T)$
denote the set of bars in $\ladder$ that $X$ has crossed during $[0,t]$. Let the set of time-$t$ {\em untouched} bar locations $\ubl_t \subseteq E(\treeg) \times [0,T)$ 
denote the set of bars $b \in E(\treeg) \times [0,T)$ neither of whose joints belongs to $X_{[0,t]}$. Then the conditional distribution of $\ladder$ is given by $\forced_t \cup \ladder_{(t,\infty)}$, where $\ladder_{(t,\infty)}$ is a random bar collection with Poisson law of intensity~$1\!\!1_{\ubl_t}$ with respect to Lebesgue measure on $E(\treeg) \times [0,T)$.
\end{lemma}
\noindent{\bf Proof.} That $\forced_t \subseteq \ladder$ is known given $X$ on $[0,t]$; similarly, if $X_{[0,t]}$ visits the joint of some bar in $\ladder$, 
that bar belongs to $\forced_t$. The time-$0$ distribution of the remaining bars, those in $\ubl_t$, is undisturbed by the data~$X_{[0,t]}$. \qed

\subsection{Useful bars}\label{secusefulbars}

We outline the strategy for proving Proposition~\ref{propone}.
As time evolves from the outset, the process $X$ will, with positive probability, jump across several bars in $\ladder$, so that $Y$ may start to move away from $\phi$. For each $t > 0$,
we will identify a subset $\good_t$ of $\forced_t$ of ``useful'' bars that, roughly speaking, act as regeneration points for the trajectory $Y:[0,t] \to V(G)$. Among other properties, these bars  have  been crossed by $X$ only once before time $t$, so that, subsequently to crossing a useful bar, the walk $Y$ 
is a descendent of the child vertex of the edge on which the bar is supported.
If the walk is to return to $\phi$, it must later pass back along each edge that supports a useful bar; however, we will choose a definition of useful bar so that, in endeavouring to return, the walk necessarily runs a positive probability of jumping out into a previously unvisited sub-tree. If the walk arrives in such a sub-tree, we will argue that there is a significant chance that it moves forward from there for a short while, thereby generating further useful bars (the number of which grows linearly with the period $T$). 
That is, to return to the root, the walk must ``undo'' each useful bar; but any attempt to do so will generate many more such bars with a uniformly positive probability. This means that the walk returns to the root only with small probability, if $T$ is chosen to be high.     

We now identify the subset $\good_t \subseteq \forced_t$. Some more notation is needed.
\begin{definition}
For any subset $A \subseteq V(G)$, let $H_A \in [0,\infty]$,  $H_A = \inf \big\{ t \geq 0: Y(t) \in A \big\}$,
denote the hitting time of $A$ by $Y$; the convention $\inf \emptyset = \infty$ is used. 
\end{definition}
Let $t > 0$ and take $B \in \forced_t$, with $B = (e,s) \in  E(G) \times [0,T)$.
Then we define $B$ to be an element of $\good_t$ if each of the following conditions is satisfied:
\begin{itemize}
\item $H_{\parent{e}} < H_{\child{e}} < t$;
\item $\big\{ s \in [0,t]: Y(s) = \parent{e} \big\} = [H_{\parent{e}},H_{\child{e}})$; 
\item $H_{\child{e}} - H_{\parent{e}} \leq \kappan$; and
\item the set $\big\{ s \in [0,t): Y(s) = \child{e} \big\}$ takes the form of an interval whose right-hand endpoint is strictly less than $t$. 
\end{itemize}
In other words, a bar $(e,s) \in \ladder$ crossed before time $t$ is useful (at time $t$) if, in its history strictly before time $t$, the walk $Y$ has made a jump and arrived at the edge $e$'s parent vertex $e^+$, 
and has then, without intervening jumps and before a duration $\kappan$ has passed, jumped to the child vertex $e^-$,
before jumping again to one of the offspring of $e^-$, without then returning to $e^-$.

\subsection{The return to a useful bar}

Each vertex in $G$ having at least two offspring, we note that, for $(e,s) \in \good_t$, each of $e^+$ and $e^-$ has an offspring, $u^+$ and $u^-$, such that
 $Y_{[0,t]} \cap \big\{ u^+, u^- \big\} = \emptyset$; indeed, there are $\degr{\parent{e}} - 2 \geq 1$ choices for $u^+$ and $\degr{\child{e}} - 2 \geq 1$ for $u^-$, because $Y$ until time $t$ has visited at most one offspring of $\parent{e}$ and at most one offspring of $\child{e}$. This fact explains how we will be able to treat the elements of $\good_t$ as obstacles for the return of $Y$ to $\phi$ after time $t$. We will argue that, conditionally on returning to $e^-$ after time $t$, there is positive probability that $Y$ arrives at either $u^+$ or $u^-$. 

We make some definitions before proving a lemma to this effect.


\begin{definition}
The time 
$t \in [0,\infty)$ 
is called a frontier time  of $X:[0,\infty) \to V(G) \times [0,T)$ if $Y(t) \not\in Y_{[0,t)}$. 
\end{definition}
\begin{definition}
For $t \geq 0$ and $A \subseteq V(G)$, let $H_{t,A} \in [t,\infty]$ be given by $H_{t,A} = \inf \big\{ s \geq t: Y(s) \in  A \big\}$. For $x \in V(G)$, we write $H_{t,x} = H_{t,\{ x\}}$.
\end{definition}

\begin{definition}
Let $t > 0$. Consider $\PR_T$ given $X:[0,t] \to V(G) \times [0,T)$ and let $(e,s) \in \good_t$. 
Also given $H_{t,\child{e}} < \infty$, we say that $X$ makes a frontier departure from $e$  if, 
after time $H_{t,\child{e}}$, at the moment of departure of $Y$ from $\{ e^+,e^- \}$, $Y$ arrives at 
an offspring of either $\parent{e}$ or $\child{e}$ that it has never visited before.
\end{definition}

\begin{lemma}\label{lemnew}
Assume~(\ref{eqatleastdegg}).
Let $t > 0$. Consider $\PR_T$ given $X:[0,t] \to V(G) \times [0,T)$; choosing $(e,s) \in \good_t$ such that $e^+ \not= \phi$ measurably with respect to $X_{[0,t]}$, condition   
further on $H_{t,\child{e}} < \infty$. Then  
$X$ makes a frontier departure from $e$ with probability at least $\frac{\degg - 1}{\degg + 1} \big( 1 - e^{-(\degg - 1)\tminusk}\big)$.
\end{lemma}

\begin{figure}
\centering\epsfig{file=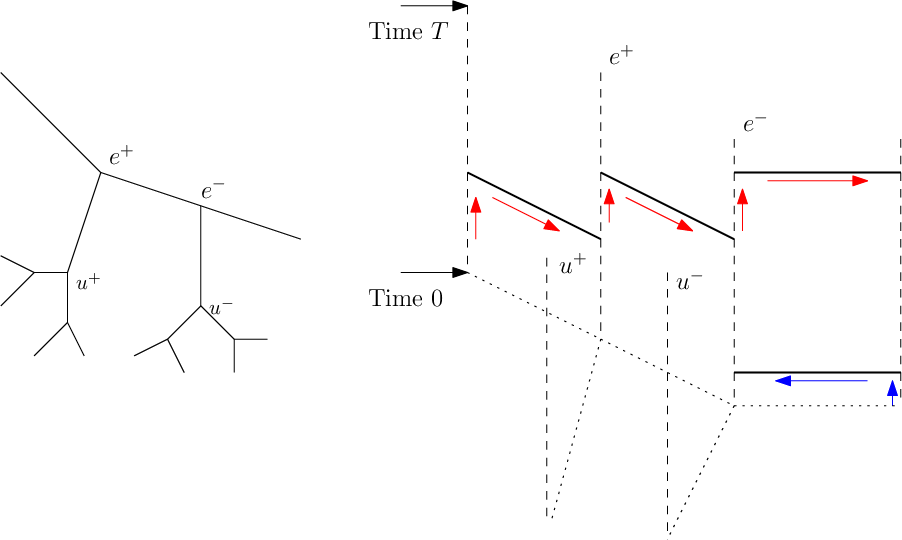, width=14cm}
\caption{The return of $X$ to the pole at the child vertex of an edge $e$ supporting a bar in $\good_t$ is depicted. In the left-hand figure, the locale of $G$ near $e$ is shown, including the beginnings of the descendent trees of $u^+$ and $u^-$; correspondingly, the dotted line segments at the base of the right-hand sketch indicate the graph structure.  The red arrows indicate the trajectory of the meander $X$ until time~$t$. The blue arrows indicate the trajectory of $X$ just prior to its return to the pole of $\child{e}$ at time~$H_{t,\child{e}}$.}\label{figreturntoedgee}
\end{figure}

\noindent{\bf Proof.}
Let $U^+$ (and $U^-$) denote the set of offspring of $\parent{e}$ (and $\child{e}$) that $Y$ has not visited during $[0,t]$.
Since the only offspring of $e^+$ that $Y$ has visited by time $t$  is $\child{e}$, and there is only one offspring of $\child{e}$ that has been so visited, we have that 
\begin{equation}\label{equplusuminus}
   \vert U^+ \vert \geq \degr{e^+} - 2 \, \quad \textrm{and} \quad \vert U^- \vert \geq \degr{e^-} - 2 \, .
\end{equation}
  
In Figure \ref{figreturntoedgee}, $u^+ \in U^+$ and $u^- \in U^-$ have been fixed for the purpose of illustration. Writing $\hnew = H_{t,\child{e}}$, we will apply Lemma \ref{lemubl} to study the conditional distribution of $X(\hnew + \cdot):[0,\infty) \to V(G) \times [0,T)$.
Note that $(\child{e},s)$ is the first joint of a bar in $\forced_\hnew$ encountered on the pole at $\child{e}$  by unit-speed cyclic upward motion  from $X(\hnew)$. Indeed, $X$ has visited the pole at $\child{e}$ before time $\hnew$ during a single interval of time, arriving there at the joint $(\child{e},s)$; now, returning to $\child{e}$ at time $\hnew$,
this point is the first joint of an element of  $\forced_\hnew$ to be located on a journey upwards from~$X(\hnew)$.

Note also that, given $X_{[0,\hnew]}$, were $X$ after time $\hnew$  to remain at $\child{e}$ until encountering the joint $(\child{e},s)$, it would jump to $(\parent{e},s)$ at the moment of this encounter. 
Let the ``jump" event $J$ occur if $X$ does indeed remain at the pole at $\child{e}$ until meeting $(\child{e},s)$; let $I \subseteq [0,T)$ denote the interval of heights 
through which $X$ passes after time $\eta$ and before encountering $(\child{e},s)$ in the case that $J$ occurs.

Note that $X_{[0,\hnew]} \cap \big( U^- \times [0,T) \big) = \emptyset$, so that $\{ (e^-,u^-) \} \times I \subseteq \ubl_{\hnew}$ for every element $u^-$ of $U^-$.
In the notation of Lemma~\ref{lemubl}, under $\PR_T$ given $X_{[0,\hnew]}$, $J^c$ occurs if and only if $\ladder_{(\eta,\infty)}$ contains a bar with a joint in $\{ e^- \} \times I$;
the fact that $\{ (e^-,u^-) \} \times I \subseteq \ubl_{\hnew}$ for every $u^- \in U^-$ and~(\ref{equplusuminus}) thus ensures that,  under $\PR_T$ given $X_{[0,\hnew]}$ and $J^c$, there is probability at least $1 - 2/\degr{e^-}$ 
that the element of $\ladder_{(\eta,\infty)}$ of lowest height among those having a joint in $\{ e^- \} \times I$ is supported on $(e^-,u^-)$ for some $u^- \in U^-$. 
 In this way, we see from Lemma~\ref{lemubl} and the strong Markov property that
\begin{equation}\label{eqescone}
 \PR_T \Big(  Y (\tau) \in U^-  \, \Big\vert \, X_{[0,\hnew]}, J^c  \Big) \geq \frac{\degr{e^-} - 2}{\degr{e^-}} \, .
\end{equation}

On the other hand, conditionally on $X_{[0,\hnew]}$ and on $J$, $X$ after time $\hnew$  leaves the pole at $\child{e}$ by crossing the bar $(e,s)$ to arrive at $(\parent{e},s)$. 
Let $\chi = H_{ \eta , V(G) \setminus \{ \child{e} \}}$ 
denote the moment of this arrival. Noting that $(e,s) \in \good_t$ and that $X_{[t,\chi)} \cap \big( \{ \parent{e} \} \times [0,T) \big) = \emptyset$, we see that 
$\big( \{ \parent{e} \} \times [0,T) \big) \cap X_{[0,\chi)}$ consists of an interval of length at most $T/2$ whose upper endpoint is $X(\chi)$, so that 
\begin{equation}\label{eqvisitchi}
(\parent{e},r \, {\rm mod} \, T) \not\in X_{[0,\chi]}
\textrm{ for all $r \in ( \chi,\chi + \tminusk )$} \, .
\end{equation} 
Moreover, $\big(  U^+  \times [0,T) \big) \cap X_{[0,\chi]} = \emptyset$.
In light of these facts, we will see that Lemma \ref{lemubl} implies that
\begin{equation}\label{eqesctwo}
 \PR_T \Big(  Y (\tau) \in U^+  \, \Big\vert \, X_{[0,\hnew]}, J  \Big) \geq  \frac{\degr{e^+} - 2}{\degr{e^+}} \Big( 1 - e^{-(\degg - 1)\tminusk}\Big) \, .
\end{equation}
Indeed,
 (\ref{eqvisitchi}) implies that unit-speed upward cyclic motion  from $(e^+,s)$
will meet no joint of a bar in $\forced_\chi$ for a duration of at least $\tminusk$. 
Whenever $u^+ \in U^+$, which is to say that $u^+$ is a neighbour of $e^+$ such that $u^+ \not\in Y_{[0,\chi)}$, $(e^+,u^+) \times [0,T) \subseteq \ubl_\chi$. 
Since any offspring of $e^+$ except for $e^-$ belongs to $U^+$, we have that  $\vert U^+ \vert \geq \degg - 1$. Applying Lemma~\ref{lemubl}, 
the conditional probability that there exists a bar in $\ladder \setminus \forced_\chi = \ladder_{(\chi,\infty)}$ with a joint at $Y(\chi)$ having a height lying in the modulo-$T$ reduction of $(\chi,\chi + \tminusk)$  is at least 
$1 - \exp\big\{ - (\degg - 1) \tminusk \big\}$; note also that, by~(\ref{equplusuminus}), 
the conditional probability, 
given the presence of such a bar, 
that the first such bar encountered by upward cyclic motion from $X_\chi$ is supported on an edge of the form $(e^+,u^+)$ for some $u^+ \in U^+$ is at least $1 - 2/\degr{e^+}$. Hence, we obtain (\ref{eqesctwo}) by applying  Lemma \ref{lemubl} (and the strong Markov property) at time $\chi$.

The lemma follows from (\ref{eqescone}) and (\ref{eqesctwo}), since $\min \big\{ \degr{e^+},\degr{e^-} \big\} \geq \degg + 1$. \hfill $\Box$

\subsection{Departing after the return}

We now extend the notion of a useful bar, making it relative to a non-zero start time. 
\begin{definition}
Let $s > 0$, and let $t > s$. Let $\forced_{s,t} \subseteq E(G) \times [0,T)$ denote the set of bars crossed by $X$ during $[s,t)$.

Let $(e,r) \in \forced_{s,t}$.
We declare that  $(e,r) \in \good_{s,t}$ if each of the following conditions is satisfied:
\begin{itemize}
\item $\parent{e}$ is a strict descendent of~$Y(s)$;
\item $d \big( Y(s), Y(H_{s,\parent{e}}) \big) > d \big( Y(s), Y(r) \big)$ for all $r \in [s,H_{s,\parent{e}})$;
\item  $\big\{ r \in [s,t]: Y(r) = \parent{e} \big\} =  [H_{s,\parent{e}},H_{s,\child{e}})$; 
\item $H_{s,\parent{e}} - H_{s,\child(e)} \leq \kappan$; and 
\item the set $\big\{ r \in [s,t): Y(r) = \child{e} \big\}$ takes the form of an interval whose right-hand endpoint is strictly less than $t$. 
\end{itemize}
\end{definition}
In fact, the set $\good_{0,t}$ may be smaller than $\good_t$, because we do not require the second of the above conditions in defining $\good_t$.
The stricter definition permits the union property recorded in the next lemma. 
\begin{definition}
Let $t > 0$. A time $s \in [0,t]$ is called a $t$-regeneration time if $\big\{ r \in [0,t]: Y(r) = Y(s) \big\}$
is an interval. If this interval has right-hand endpoint strictly less than $t$, then the first jump made by $Y$ after time $s$ is in the direction away from the root. 
\end{definition}
\begin{lemma}\label{lemdisjoint}
Let $t > s > 0$.
 Let $X:[0,\infty) \to V(G) \times [0,T)$ be such that $s$ is a $t$-regeneration time. Then  $\good_s$ and $\good_{s,t} $ are disjoint subsets of $\good_t$.  
\end{lemma}
\noindent{\bf Proof.} Let $(e,r) \in \good_s$, (with $e \in E(G)$ and $r \in [0,T)$). Note that $Y_{[s,t]}$ lies in the descendent tree of $Y(s)$, a vertex which is itself a strict descendent of~$e^-$. 
Hence, $\good_s \subseteq \good_t$.
The other inclusion is similarly established. No vertex associated to a joint of a bar in $\good_s$ is a  descendent of $Y(s)$, while every vertex associated to a joint of a bar in $\good_{s,t}$ is a strict descendent of $Y(s)$. This ensures the disjointness of the two sets. \hfill $\Box$

\begin{lemma}\label{lemrapidadvanceone}
Assume~(\ref{eqatleastdegg}).
Given $\epsilon > 0$, there exists  $T_0 > 0$ such that, for $T \geq T_0$,  
the $\PR_T$-probability that 
\begin{equation}\label{gbound}
\vert \good_{0,T} \vert \geq 
 \bigg(  \frac{\degg ( \degg - 1 )^2}{(\degg + 1)^2}    \Big( 1 - e^{- ( \degg + 1 ) \kappan} \Big) 
- \epsilon \bigg) T
\end{equation}
is at least $1 - \epsilon$.
\end{lemma}
\noindent{\bf Proof.} 
Let $\tdegg$ denote the rooted regular tree each of whose vertices has $\degg$ offspring.
We first prove the lemma when $G = \tdegg$.

 Let $\beta > 1$.
Let $Z:\N \to \N$, $Z(0) = 0$, denote nearest-neighbour random walk with bias $\beta$ to the right and with reflection at zero.
This is the Markov chain with transition probabilities $p_{n,m} = \delta_{m,n+1} \frac{\beta}{\beta + 1} + \delta_{m,n-1} \frac{1}{\beta + 1}$ for $n \geq 1$ and $p_{0,m} = \delta_{m,1}$. 
We call  $n \in \Z$ a renewal point for $Z$  (and write $n \in \zreg$)
if $m \in \N$ and $Z(m) = Z(n)$ implies that $m =n$.  We call $n$ a strong renewal point for $Z$ (and write $n \in \zstrreg$) if $\{ n,n+1\} \subseteq \zreg$. Note that the conditional distribution given $Z:[0,n] \to \Z$
and $n \in \zstrreg$ of $Z(n + \cdot) - Z(n)$ is given by $Z$ conditioned to make three rightward steps and then to remain at values of at least two. This conditional distribution being independent of  $Z:[0,n] \to \Z$ given 
$n \in \zstrreg$, we see that the strong renewal points form a renewal sequence (in the sense that the differences of consecutive terms are independent and have a common law).
It is easy to confirm that, for each $n \in \N^+$, $\mathbb{P} (n \in \zstrreg) = \frac{\beta (\beta - 1)}{(\beta + 1)^2}$. Hence, the renewal theorem implies that
\begin{equation}\label{eqzstrbeta}
   n^{-1} \Big\vert \zstrreg \cap \big\{ 1,\ldots, n \big\} \Big\vert \to \frac{\beta ( \beta - 1 )}{(\beta + 1)^2}, \qquad \textrm{almost surely.}
\end{equation}

Let $W:[0,\infty) \to V(G)$ denote continuous-time random walk on $G$ departing from $\phi$ (whose jumps are given by exponential rate-one clocks on the edges of $G$).
Write $M:[0,\infty) \to \N$, $M(s) = d \big( \phi , W(s) \big)$, where recall that $d(\cdot,\cdot)$
denotes graphical distance on $G$. Let $J:\N \to \N$, $J(0) = 0$, denote the jump chain of $M$, that records in discrete-time the successive states visited by $M$. Let $D:\N \to (0,\infty)$,
where $D(0)$ is the time for $W$ to make its first transition, and $D(n)$ for $n \in \N^+$
is the length of time that $W$ spends at its new location after its $n$-th transition.

Taking $\beta = d$, note that $J:\N \to \N$ and $Z_\beta:\N \to \N$ are equal in law.

Note that $Y:[0,T) \to V(G)$ has the distribution of $W:[0,T) \to V(G)$.
Thus, we wish to argue that~(\ref{gbound}) holds for the process  $W:[0,T) \to V(G)$.  
The process~$W$ making transitions at rate at least $\degg + 1$ except when at the root, the law of large numbers implies that, for any $\epsilon > 0$, there exists an almost surely finite random variable $T_0$ such that   
 $W:[0,T) \to V(G)$  makes at least $T(\degg + 1)(1 -\epsilon)$ transitions if $T \geq T_0$. Clearly then, for $\epsilon > 0$, we have that $W(T) \geq T (\degg - 1) (1 - \epsilon)$ if $T \geq T_0$ (where the law of the almost surely finite $T_0$ may have changed).
 Every strong renewal point $j$ for $J$ for which $0 < j \leq T(\degg - 1) (1-\epsilon)$ corresponds to an element of  $\good_{0,T}$, provided that 
$D(j) \leq \kappan$.  For any given $k \in \N^+$, conditionally on a choice of $J:\N \to \N$ such that $J(k) > 0$, and on the values $\big\{ D(i): i \not= k \big\}$, the conditional probability that $D(k) \leq \kappan$ is at least $1 - e^{- (\degg + 1) \kappan}$.

Recalling (\ref{eqzstrbeta}), we see that, for any $\epsilon > 0$, there exists a deterministic $T_0 > 0$ such that  $T \geq T_0$ implies that
$$
  \big\vert \good_{0,T} \big\vert \geq \bigg(  \frac{\degg ( \degg - 1 )}{(\degg + 1)^2} - \epsilon \bigg) T  (\degg - 1)   \big( 1 -\epsilon \big) \Big( 1 - e^{- (\degg +1) \kappan} \Big)
$$
with probability at least $1 - \epsilon$.
This yields the statement of the lemma in the case that $G = \tdegg$.

The general case may be reduced to the special one. Assume now that $G$ satisfies~(\ref{eqatleastdegg}).
It is straightforward to construct a coupling $\coupling$
of the continuous-time random walks $W^{\tdegg}$ and $W^G$ with the property that
at any moment of time at which $W^G$ makes a transition towards the root, $W^{\tdegg}$ is either at the root or makes a transition towards the root, while at any moment of time at which $W^{\tdegg}$ makes a transition away from the root, so does $W^G$. We omit the details of this standard construction. Under $\coupling$, let $T_1$ denote the supremum of times at which $W^{\tdegg}$ is at the root, and let $m$ denote the maximal distance from the root attained by $W^G$ before $T_1$. Let $T_2 = \inf \big\{ t \geq T_1: d \big( \phi, W^G \big) = m \big\}$. 
It is easily verified that, 
if $T > T_2$ and $s \in (T_2,T)$ is a moment 
at which a bar in $\good_{0,T}^{\tdegg}$ is crossed by $W^{\tdegg}$,
then $s$ is also a moment at which
a bar in $\good_{0,T}^G$ is crossed by $W^G$. 
Hence, whenever $T > T_2$, $\big\vert \good_{0,T}^G \big\vert \geq \big\vert \good_{0,T}^{\tdegg} \big\vert - \big\vert \good_{0,T_2}^{\tdegg} \big\vert$ under $\coupling$.  The random variable $T_2$ being finite $\coupling$-almost surely, we infer the statement of the lemma for the graph $G$ from this statement for $\tdegg$. \hfill $\Box$

\begin{lemma}\label{lemrapidadvance}
Assume~(\ref{eqatleastdegg}).
Given $\epsilon > 0$, there exists  $T_0 > 0$ such that the following holds.
Let $T \geq T_0$. Fix $t > 0$. Consider $\PR_T$ given  $X:[0,t] \to V(G) \times [0,T)$
such that $t$ is a frontier time. 
Then the conditional probability that 
\begin{equation}\label{eqgoodbd}
\vert \good_{t,t+T} \vert \geq
 \bigg(  \frac{\degg ( \degg - 1 )^2}{(\degg + 1)^2}    \big( 1 - e^{- (\degg + 1 ) \kappan} \big) 
- \epsilon \bigg) T
\end{equation}
and that $t$ is a $(t+T)$-regeneration time 
is at least $1/4$.
\end{lemma}
\noindent{\bf Remark.} 
Define 
$$
\cbeta = 
\frac{\degg ( \degg - 1 )^2}{2(\degg + 1)^2}    \Big( 1 - e^{- (\degg + 1) \kappan} \Big).
$$ 
Let $X:[0,\infty) \to V(G) \times [0,T)$ 
be such that 
$t \in [0,\infty)$ 
is a frontier time of~$X$.
We say that $X$ makes a rapid advance from $t$ if $Y(s)$ is a descendent of $Y(t)$ for all $s \in [t,t+T]$
and if  
$\vert \good_{t,t+T} \vert \geq \cbeta T$. 
In these terms, Lemma~\ref{lemrapidadvance} implies that there exists $T_0$ such that, if $T \geq T_0$, then, 
under $\PR_T$  given $X:[0,t] \to V(G) \times [0,T)$ such that $t$ is a frontier time, 
the conditional probability that $X$ makes a rapid advance from $t$ is at least $1/4$.

\noindent{\bf Proof of Lemma \ref{lemrapidadvance}.}
Let $D$ denote the descendent tree of $Y(t)$, (with root $Y(t)$). 
In this proof, we use the term $G$-process to refer to the conditional distribution of $X(t + \cdot)$
given  data $X:[0,t] \to V(G) \times [0,T)$ as specified in the statement of the lemma. By the $D$-process, we mean cyclic-time random meander in $D$ from $X(t)$. In this way, each process is defined on $[0,\infty)$.
Note that the $D$- and $G$-processes may be coupled 
 by using the same collection of elements in $\ladder$ supported on edges in $D$. 
The two processes then coincide during $[0,T]$
provided that the $G$-process departs from the pole at $Y(t)$ by jumping to the pole indexed by an offspring of $Y(t)$, and does not return to the pole at $Y(t)$ during $[0,T]$. 
The probability that this happens is at least 
$\frac{\degg - 1}{\degg + 1}$ by Lemma \ref{lemubl} and a basic hitting estimate.
By Lemma \ref{lemrapidadvanceone}, the random variable $\good_{t,t+T}$ under the $D$-process satisfies (\ref{eqgoodbd}) with probability arbitrarily close to one (provided that $T$ is chosen to be high enough).
If the coupling works, then the bound applies to the $G$-process as well. The bound $\degg \geq 2$ gives the statement of the lemma (where in fact $1/4$ might be replaced by any value less than $1/3$).    \hfill $\Box$


\begin{definition}\label{defgoodreturn}
Let $t > 0$ and let the bar $(e,r)$ denote any given element of~$\good_t$. 
Write $e^c = V(G) \setminus \big\{ \parent{e},\child{e} \big\}$.
We say that $X$ makes a return to $e$ if $H_{t,\child{e}} < \infty$. If $X$ makes a return to $e$, we say that the return is good if
\begin{enumerate} 
\item $X$ makes a frontier departure from $e$, and then 
\item $X$ makes a rapid advance from the frontier time $H_{H_{t,\child{e}},e^c}$.
\end{enumerate}
\end{definition}
Define
$$
\ctbeta =   \frac{\degg - 1}{4(\degg + 1)}  \big( 1 - e^{-(\degg - 1)\tminusk}\big) \, .
$$
\begin{lemma}\label{lemsum}
Assume~(\ref{eqatleastdegg}).
There exists $T_0 > 0$ such that the following holds. Let $T \geq T_0$.
For any $t > 0$ and any $(e,r) \in E(G) \times [0,T)$ with $e^+ \not= \phi$, under $\PR_T$ given
  $(e,r) \in \good_t$ and $H_{t,\child{e}} < \infty$,
the probability that the return of $X$ to $e$ is good is at least~$\ctbeta$.
\end{lemma}
\noindent{\bf Proof.}
 This is implied by  Lemma \ref{lemnew}
and the remark that follows Lemma~\ref{lemrapidadvance}, as well as by Lemma \ref{lemubl}. \hfill $\Box$

\subsection{Damage limitation after a bad return}

We need a lemma that controls the damage done by a return to a useful bar which turns out not to be good.
\begin{lemma}\label{lemretnotgood}
 Let $t > 0$. Let $e$ denote the element of $\good_t$ that is crossed last by $X:[0,t] \to V(G) \times [0,T)$.
Let $p(e^+)$ denote the parent of $e^+$. Then, conditionally on $e^+ \not= \phi$ and $H_{t,p(e^+)} < \infty$, we have that,  
almost surely,
 $\good_t \setminus \good_{H_{t,p(e^+)}}$  contains at most two elements.
\end{lemma}
\noindent{\bf Proof.} We write $\overline\good_t$ for the set of edges that support a bar in $\good_t$.
Note that, for each $t > 0$, these two sets are in one-to-one correspondence, since no two elements in $\good_t$ are supported on the same edge.
 Hence, it suffices to derive the statement of the lemma with $\good$ replaced by $\overline\good$.

Note that the elements of $\overline\good_t$, enumerated $\big( e_1,\ldots, e_k \big)$
in the order in which the constituent edges are crossed by $Y:[0,t] \to V(G) \times [0,T)$,
have the property that, in the list  $\big( e_1^-,\ldots, e_k^- \big)$,
each entry is a descendent of each of its precursors. By definition, $e = e_k$.
For $1 \leq i \leq k-1$, note that $\inf \big\{ s \geq t: e_i \not\in \overline\good_s \big\} = H_{t,e_i^-}$. Note also that if $H_{t,e_i^-} < \infty$ then $e_i \in \overline\good_{H_{t,e_i^-}}$, 
because the fourth requirement in the definition of $\{ \good_t:t\geq 0 \}$ in Subsection~\ref{secusefulbars} is chosen so that a given bar leaves this process not at, but only momentarily after, a return by $Y$ to the child vertex of the edge supporting this bar.  
Thus, $\overline\good_t \setminus \overline\good_{H_{t,p(e^+)}}$ may contain no edges other than $e$ and $\big(p(e^+),e^+\big)$.  \hfill $\qed$

\subsection{Establishing the main results}

$\empty$

\noindent{\\ {\bf Proof of Proposition~\ref{propone}.}}
We fix $T_0 > 0$ high enough to satisfy the hypotheses of each of the preceding lemmas.
We now form the process $X:[0,\infty) \to V(G) \times [0,T)$ iteratively. 
We will construct an increasing sequence $\big\{ \tau_i : i \in \N^+ \big\}$ 
of times  at which the present number of useful bars is gauged.

At the first step, provided that the positive probability event that 
$\vert \good_t \vert \geq 2$ for some $t > 0$ occurs, we set  
$\tau_1 = \inf \big\{ t \geq 0 : \vert \good_t \vert \geq 2 \big\}$. (We wait for two elements to appear in $\good_t$ because if there is only one, it may be supported on an edge incident to $\phi$, and we have not set up the tools to handle this case.) Otherwise, we set $\tau_j = -\infty$ for each $j \in \N^+$, in a formal device indicating that our effort to determine that $X$ does not have a periodic trajectory has failed. 
(In the case that $\tau_1 \not= - \infty$, note that $\vert \good_{\tau_1} \vert = 2$, because the convention that $X$ be right-continuous means that $X$ crosses a bar at time $\tau_1$, and is no longer at the child vertex of the second bar to become useful, permitting this bar to join $\good_t$ at time $t = \tau_1$.)

Let $k \in \N^+$. Suppose that $0 < \tau_k < \infty$. (As will be apparent, we set $\tau_k = \infty$ in the case that it becomes evident at the $k$-th stage of the construction that $X$ does not have a periodic trajectory. For definiteness, if $\tau_k$ is set equal to $\infty$ in the subsequent definition, then we automatically also set  $\tau_l = \infty$ for all $l > k$.)

If  $\vert \mathcal{U}_{\tau_{k}} \vert \leq 1$, set $\tau_l = - \infty$ for all $l \geq k+1$. 

Otherwise, let $(e_k,t_k)$ denote the bar in $\good_{\tau_k}$ that is the last to be crossed by $X$ 
before time $\tau_k$.
Let $\chi_k = H_{\tau_k,e_k^-}$.

If $\chi_k = \infty$, then set $\tau_{k+1}=\infty$.

If $\chi_k < \infty$ and the return of $X$ to $e_k$ is good, 
recalling that $e_k^c$ denotes $V(G) \setminus \{ e_k^+ , e_k^- \}$, we take $\tau_{k+1}= H_{\chi_k,e_k^c} + T$. 

If $\chi_k < \infty$ and the return of $X$ to $e_k$ is not good, we take $\tau_{k+1} = H_{\chi_k,p(\parent{e})}$, 
where recall that $p(\parent{e})$ denotes the parent of $\parent{e}$. Note that  $H_{\chi_k,p(\parent{e})}$  may be infinite.

For $k \in \N^+$, set $u_k = \vert \mathcal{U}_{\tau_k} \vert$. As a convention, we take $u_k = 0$ if
$\tau_k = -\infty$ and $u_k = \infty$ if $\tau_k = \infty$.

For $t \in [0,\infty)$, let $\sigma_t$ denote the sigma-algebra generated by $\big\{ X_s: 0 \leq s \leq t \big\}$.
For $k \in \N^+$, write $\sigma'_k = \sigma_{\tau_k}$, where, in a standard definition,
$\sigma_{\tau_k} = \big\{ A \subseteq \Omega: A \cap \big\{ \tau_k \leq t \big\} \in \sigma_t \, \textrm{for each $t  > 0$} \big\}$.

We now define three $\sigma'_k$-measurable random variables, $p_k$, $q_k$ and $r_k$.
To define each of them, consider $\PR_T$ given  $\big\{ X_t: 0 \leq t \leq \tau_k \big\}$.
Then $p_k$ is set  equal to the conditional probability that $\chi_k < \infty$.
Let $q_k$ denote the conditional probability,  given further that $\chi_k < \infty$, that the return of $X$ to $e_k$ is good. Let $r_k$ denote the conditional probability,  given that $\chi_k < \infty$ and that the return of $X$ to $e_k$ is not good, that $H_{\chi_k,p(\parent{e})} < \infty$.

Note that $u_k$ is $\sigma'_k$-measurable. Note also that, by the definition of a good return, and Lemmas~\ref{lemdisjoint} and~\ref{lemretnotgood}, we have that, $\sigma'_k \big( \cdot \big\vert u_k > 1 \big)$-almost surely, the conditional distribution of $u_{k+1} - u_k$ given $\big\{ X_t:0\leq t \leq \tau_k \big\}$ stochastically dominates the
law $\infty \delta_{1 - p_k} + \cbeta T \delta_{p_k q_k} + \infty \delta_{p_k(1 - q_k)(1-r_k)}  -2   \delta_{p_k(1 - q_k)r_k}$. The latter distribution is parametrized by $(r_k,p_k,q_k)$, and stochastically dominates the one obtained by replacing the values of each of $p_k$ and $r_k$ by $1$.
Moreover, by Lemma~\ref{lemsum}, $q_k \geq \ctbeta$, $\sigma'_k \big( \cdot \big\vert u_k > 1 \big)$-almost surely.
To summarise these  deductions,  $\sigma'_k  \big( \cdot \big\vert u_k > 1 \big)$-almost surely, the conditional distribution of $u_{k+1} - u_k$ given $\big\{ X_t:0\leq t \leq \tau_k \big\}$ stochastically dominates the law 
$\cbeta T \delta_{\ctbeta} - 2 \delta_{1 - \ctbeta}$.

The data $\big\{ u_1,\ldots,u_k \big\}$ being $\sigma_k'$-measurable, we infer that, given such data for which $u_k > 1$, the conditional distribution of $u_{k+1} - u_k$
also stochastically dominates the law 
$\cbeta T \delta_{\ctbeta} - 2 \delta_{1 - \ctbeta}$.
Let $Q:\N^+ \to \R$ denote the random walk on $\R$ whose increments are independent and have the law
$\cbeta T \delta_{\ctbeta} - 2 \delta_{1 - \ctbeta}$, with initial condition $Q(1) = 2$.
Let $\rho \in \N$ denote the first time at which $Q$ is at most one, and define $Q_*:\N^+ \to \R$ by
\begin{equation*}
Q_*(i)  =  \begin{cases}  Q(i) &  \text{if } i \leq \rho \, ,  \\
            0 & \text{if } i > \rho \, ,
             \end{cases}
 \end{equation*}
for each $i \in \N^+$.
 
We find that, conditionally on $\tau_1 < \infty$,  $\big\{ u_i: i \in \N^+ \big\}$
stochastically dominates $\big\{ Q_*(j): j \in \N^+ \big\}$.
Hence, we find that the probability that $u_i \to \infty$ as $i \to \infty$
is at least the probability that $\big\{ Q(i): i \in \N \big\}$ is a sequence of terms all of which exceed one and which tend to infinity.
By the law of large numbers,
this occurrence has  positive probability provided that 
\begin{equation}\label{eqcdtcb}
  \ctbeta \cbeta T - 2 \big( 1 - \ctbeta \big) > 0 \, .
\end{equation}
Note that the left-hand side is non-decreasing and tends to infinity in the limit of high $T$; as such, we may adjust the value of $T_0 > 0$, if necessary, so that, for $T > T_0$, the condition (\ref{eqcdtcb}) is satisfied.

Clearly, if $X:[0,\infty) \to V(G) \times [0,T)$ has a periodic orbit, then $\tau_k$ eventually assumes the value $-\infty$. 
Hence, provided that $T> T_0$, with positive probability, $X$ does not have a periodic orbit,  so that $\phi \not\in Y(t,\infty)$ for sufficiently high $t$. This completes the proof of Proposition~\ref{propone}.  \hfill $\Box$

\vspace{1mm}

\noindent{\bf Proof of Theorem~\ref{thmone}.} 
By Proposition~\ref{propone}, there is positive probability that $X = X_{(\phi,0)}$ has an aperiodic orbit, in which case, members of the semi-infinite sequence $\big\{ Y(kT): k \in \N \big\}$
form the consecutive elements of part of some infinite cycle in the associated random stirring model with parameter~$T$.

Should $X$ have a periodic orbit, we may search in successive generations  away from the root for an edge $e$ such that no bar in $\ladder$ is supported on $e$. The conditional distribution of $X_{(e^-,0)}$ is then given by cyclic-time random meander in the descendent tree of $e^-$. Proposition~\ref{propone} being applicable to this tree, there is a further uniformly positive probability that $X_{(e^-,0)}$ has an aperiodic orbit. This procedure may continue until a meander with such an orbit is located. \qed

\appendix
\section{The proof of Theorem~\ref{thmtwo}}
This appendix is devoted to the proof of Theorem~\ref{thmtwo}. We begin by stating quantitative analogues to the results en route to Proposition~\ref{propone}, and then explain the changes needed to the arguments to obtain these analogues.
\subsection{Quantitative counterparts}
The first two lemmas are analogues of Lemmas \ref{lemrapidadvanceone} and \ref{lemrapidadvance}.
\begin{lemma}\label{lemappone}
Assume~(\ref{eqatleastdegg}).
If $\degg \geq 55$ and $T \geq 101 \degg^{-1}$, then
$$
\PR_T \Big(  \big\vert \good_{0,T} \big\vert \geq  \tfrac{T \degg}{9} \Big) \geq   0.649 \, .
$$
\end{lemma}
\begin{lemma}\label{lemapptwo}
Assume~(\ref{eqatleastdegg}).
Suppose that $\degg \geq 55$ and that  $T \geq 101 \degg^{-1}$. Fix $t > 0$.
Consider $\PR_T$ given  $X:[0,t] \to V(G) \times [0,T)$
such that $t$ is a frontier time. 
Then the conditional probability that 
$$
  \vert \good_{t,t+T} \vert \geq \tfrac{T\degg}{9}
$$
and that $t$ is a $(t+T)$-regeneration time 
is at least $0.626$.
\end{lemma}
We now state an analogue of the remark that follows Lemma \ref{lemrapidadvance}; we redefine the notion of ``rapid advance'' in doing so. 

\medskip

\noindent{\bf Remark.}  
Set $\cbetastar = \degg/9$.
Fix $t > 0$. Let $X:[0,\infty) \to V(G) \times [0,T)$ 
be such that 
$t \in [0,\infty)$ 
is a frontier time of $X$.
We say that $X$ makes a rapid advance from $t$ if $Y(s)$ is a descendent of $Y(t)$ for all $s \in [t,t+T]$
and if  
$\vert \good_{t,t+T} \vert \geq \cbetastar T$.
In these terms, Lemma~\ref{lemapptwo} implies that, for 
 $\degg \geq 55$ and $T \geq 101 \degg^{-1}$,  
the conditional probability given such data $X:[0,t] \to V(G) \times [0,T)$ that $X$ makes a rapid advance from $t$ is at least $0.626$.

The next result is an analogue of Lemma \ref{lemsum}.
\begin{lemma}\label{lemappthree}
Assume~(\ref{eqatleastdegg}).
Suppose that $\degg \geq 55$ and that  $T \geq 101 \degg^{-1}$. Set  $\ctbetastar =  0.626 \frac{\degg - 1}{\degg + 1}  \big( 1 - e^{-(\degg - 1)\tminusk}\big)$.
For any $t > 0$ and any $(e,r) \in E(G) \times [0,T)$ with $e^+ \not= \phi$, under $\PR_T$ given
  $(e,r) \in \good_t$ and $H_{t,\child{e}} < \infty$,
the probability that the return of $X$ to $e$ is good is at least~$\ctbetastar$.
\end{lemma}
\subsection{Proofs}

\noindent{\\ \bf Proof of Lemma \ref{lemappone}.} This is a matter of quantifying the numerous assertions made in the proof of Lemma \ref{lemrapidadvanceone}. As before, we work firstly with the tree $\tdegg$. 

Let $B_n$ denote the number of $i \in \{1,\ldots,n\}$
such that $J(i + 1) = J(i) - 1$.
Let $r \in (0,1)$ and $a > 0$ be parameters whose values will be specified later. For $t \geq 0$, 
let $\ell(t)$ denote the number of transitions made by $W$ 
during $[0,t]$. We claim that
\begin{eqnarray}\label{eqellr}
  & & \Big\{  \ell(T) \geq r T \degg + 2  a r T \tfrac{\degg}{\degg + 1} \Big\} \cap 
   \Big\{ B_{r T  \degg + 2 a r   T \degg (\degg + 1)^{-1}    } < a r T \tfrac{\degg}{\degg + 1} \Big\} \label{eqellinc} \\
  & & \cap \, \Big\{ \lfloor rT \degg \rfloor \in {\rm SRG}(J) \Big\}
 \, \subseteq \,  \Big\{  \big\vert {\rm SRG}(Z) \cap [0,T] \big\vert \geq r T \degg \big( 1 - \tfrac{3a}{\degg + 1} \Big) \Big\} \, . \nonumber 
\end{eqnarray}
Indeed, if the event on the left-hand side occurs, 
then $\sigma : = \inf \big\{ m \in \N: J(m) \geq  r T \degg \big\}$ satisfies $\sigma \leq r T  \degg + 2 a r   T \degg (\degg + 1)^{-1}$.
The process $J$ makes at most  $a r T  (\degg + 1)^{-1}$ backward steps during $[0,\sigma]$; each such step 
is responsible for at most three points being removed from ${\rm SRG}(J)$. Each element of ${\rm SRG}(J) \cap \{
1,\ldots, \lfloor r T \degg \rfloor \}$ corresponds to an element of ${\rm SRG}(W)$ of value less than $T$, since $W$ reaches $r T \degg$ before time $T$. This establishes (\ref{eqellinc}).

Note that $\ell(T)$ stochastically dominates a Poisson random variable $X$ of parameter $T \degg$. 
Hence, 
\begin{eqnarray}
  & & \mathbb{P} \Big(  \ell(T) <  r T \degg + 2  a r T \tfrac{\degg}{\degg + 1}  \Big)  \leq   \mathbb{P}  \bigg(  X <   r T \degg \Big( 1 + \tfrac{2a}{\degg+1}  \Big) \bigg) \label{eqworkone} \\
   & \leq &  \mathbb{P}  \bigg( \big\vert  X - \mathbb{E} X \big\vert \geq   T \degg \Big( 1 - r  - \tfrac{2ar}{\degg+1}  \Big) \bigg) 
   \leq   \Big(   1 - r  - \tfrac{2ar}{\degg+1} \Big)^{-2} (T \degg)^{-1} \, . \nonumber
\end{eqnarray}
Note that  $B_{r T  \degg + 2 a r   T \degg (\degg + 1)^{-1}    }$ is stochastically dominated by the binomial distribution with parameters $r T  \degg + 2 a r   T \tfrac{\degg}{\degg + 1}$ and $(\degg + 1)^{-1}$. Let $Y$ have the latter distribution. Thus,
\begin{eqnarray}
 & &  \mathbb{P} \Big(   B_{r T  \degg + 2 a r   T \degg(\degg + 1)^{-1}  } \geq a r T \tfrac{\degg}{\degg + 1}  \Big) 
  \leq  \mathbb{P} \Big(  Y   \geq a r T \tfrac{\degg}{\degg + 1}  \Big) \nonumber \\
 & \leq & \mathbb{P} \Big(  \big\vert Y - \mathbb{E}Y \big\vert   
 \geq a r T \tfrac{\degg}{\degg + 1} - \big( r T \degg + 2 a r   T \tfrac{\degg}{\degg + 1} \big) (\degg + 1)^{-1}  \Big) \nonumber \\
& \leq & \mathbb{P} \Big(  \big\vert Y - \mathbb{E}Y \big\vert   
 \geq (a - 2) r T \tfrac{\degg}{\degg + 1}   \Big) \nonumber \\
  & \leq & (a - 2)^{-2}  \Big( \tfrac{rT \degg}{\degg + 1} \Big)^{-2} {\rm Var} (Y) \leq \tfrac{4}{3} (rT)^{-1} (a-2)^{-2} \, .
   \label{eqworktwo}
   \end{eqnarray}
We imposed $2a \leq \degg + 1$ in order to obtain the third inequality and then further imposed $6a \leq \degg + 1$ to bound ${\rm Var}(Y)$. 

We now apply to (\ref{eqellinc}) the formula $\mathbb{P} \big( \lfloor rT \degg \rfloor \in {\rm SRG}(J) \big) = \tfrac{\degg(\degg - 1)}{(\degg + 1)^2}$ as well as the inequalities (\ref{eqworkone}) and (\ref{eqworktwo}) and so obtain
\begin{eqnarray}
 & & \mathbb{P}   \bigg( 
 \big\vert {\rm SRG}(Z) \cap [0,T] \big\vert \geq r T \degg \Big( 1 - \tfrac{3a}{\degg + 1} \Big) \bigg) \nonumber \\
 & \geq & \frac{\degg(\degg - 1)}{(\degg + 1)^2} - 
  \Big(   1 - r  - \tfrac{2a}{\degg+1} \Big)^{-2} (T \degg)^{-1}  - 
  \tfrac{4}{3} (rT)^{-1} (a-2)^{-2} \, . \label{eqwork}
\end{eqnarray}
Any given positive element $z$ of $\zstrreg \cap [0,T]$ belongs to $\good_{0,T}$ if and only if $D(j) \leq T/2$,
where $z$ is the $j$-th element of $Z$. Note that $\big\{ D(i): i \in \N \big\}$ forms an independent and identically distributed sequence, which is independent of $Z$, and note also that 
$\mathbb{P} \big( D(j) \leq T/2 \big) \geq 1 - e^{ - \degg T/2 }$ for each $j \in \N$.

Let $R$ have the binomial distribution of parameters $r T \degg \big(  1 - \tfrac{3a}{\degg +1} \big)$
and $1 - e^{  - \degg T/2 }$. We find then that
\begin{eqnarray}
  & & \mathbb{P} \Big(  \big\vert \good_{0,t} \big\vert + 1 < s r T \degg \big(  1 - \tfrac{3a}{\degg +1} \big) \big( 1 -  e^{  - \degg T/2 } \big) \nonumber \\
  & & \qquad \qquad \qquad \qquad  \Big\vert \, \big\vert {\rm SRG}(Z) \cap [0,T] \big\vert \geq r T \degg \big( 1 - \tfrac{3a}{\degg + 1} \big) \Big) \nonumber \\
  & \leq & \mathbb{P} \Big( R <  s r T \degg \big(  1 - \tfrac{3a}{\degg +1} \big) \big( 1 -  e^{  - \degg T/2 } \big) \Big) \nonumber \\
 & \leq & (1 - s)^{-2} \big( r T \degg \big)^{-2}   \big( 1 - \tfrac{3a}{\degg +1} \big)^{-2}
   \big( 1 -  e^{  - \degg T/2 } \big)^{-2} \, {\rm Var}(R) \nonumber \\
 & \leq &  (1 - s)^{-2} \big( r T \degg \big)^{-1}   \big( 1 - \tfrac{3a}{\degg +1} \big)^{-1}
   \big( 1 -  e^{  - \degg T/2 } \big)^{-1}  e^{  - \degg T/2 } \, .  \label{eqworknew} 
\end{eqnarray}
By (\ref{eqwork}) and (\ref{eqworknew}), 
\begin{eqnarray}
  & & \mathbb{P} \bigg(  \big\vert \good_{0,T} \big\vert \geq s r T \degg \big(  1 - \tfrac{3a}{\degg +1} \big) \big( 1 -  e^{  - \degg T/2 } \big) \, - 1 \bigg) \nonumber \\
  & \geq &   \tfrac{\degg(\degg - 1)}{(\degg + 1)^2}   
  -  \Big( 1 - r  - \tfrac{2ar}{\degg + 1} \Big)^{-2} (T\degg)^{-1} - \tfrac{4}{3} \big( r T  \big)^{-1} (a - 2)^{-1} \nonumber \\
   &  &
 - \, (1 - s)^{-2} \big( r T \degg \big)^{-1}   \big( 1 - \tfrac{3a}{\degg +1} \big)^{-1}
   \big( 1 -  e^{  - \degg T/2 } \big)^{-1}  e^{  - \degg T/2 } \, . \nonumber
\end{eqnarray}
Set $a = \tfrac{\degg + 1}{6}$, $r = s =\tfrac{1}{2}$ and write $t = T \degg$. Assuming that $\degg \geq 55$ and $t \geq 10 \log 2$, the right-hand side is at least 
$$
 1 - \tfrac{3}{\degg + 1} - \Big( 9 + 20 + 16 \cdot 2^{-5} \big(1 - 2^{-5} \big)^{-1} \Big) t^{-1} \geq  1 - \tfrac{3}{\degg + 1} -  30 t^{-1} \, .
$$
Thus, if $\degg \geq 55$ and $t \geq 10 \log 2$,
$$
\mathbb{P} \Big(  \big\vert \good_{0,T} \big\vert \geq  \tfrac{t}{8}(1 - 2^{-5}) - 1 \Big) \geq   1 - \frac{3}{\degg + 1} -  30 t^{-1} \, .
$$
From this, Lemma \ref{lemappone} follows in the case that $G = \tdegg$.

To establish the lemma for a general graph $G$, we employ the coupling $\coupling$
used in the proof of Lemma \ref{lemrapidadvanceone} to reduce the case of a general case to that of $\tdegg$.
Recall the quantities $m \in \N$, $T_1 > 0$ and $T_2 \geq T_1$ defined in the earlier argument. Consider the event
\begin{equation}\label{eqev}
  d \big( \phi, W^G (T_1) \big) = m \, ; 
\end{equation}
that is, the event that, at the moment of final departure of $W^{\tdegg}$ from the root, $W^G$
is at least as far from the root as it has ever been. Under this circumstance, we have that $T_2 = T_1$, so that, if $T > T_1$, then $\big\vert \good_{0,T}^G \big\vert \geq \big\vert \good_{0,T}^{\tdegg} \big\vert - \big\vert \good_{0,T_1}^{\tdegg} \big\vert$. However, $\good_{0,T_1}^{\tdegg} = \emptyset$, so that 
$\big\vert \good_{0,T}^G \big\vert \geq \big\vert \good_{0,T}^{\tdegg} \big\vert$ if (\ref{eqev}) holds.
In light of (\ref{eqev}), the statement of the lemma reduces to verifying that the coupling $\coupling$ may be constructed so that
\begin{equation}\label{eqcoup}
 \coupling \Big( d \big( \phi, W^G (T_1) \big) = m \Big) \geq 1 - 2\degg^{-1} \, .
\end{equation}
To show (\ref{eqcoup}), let $T_3 \in (0,T_1]$ be given by 
$T_3 = \inf \big\{ t > 0: W^{\tdegg} \not= \phi  \big\}$.
We may construct $\coupling$ so that, for each $t \geq 0$, the conditional distribution of  $W^G$ during $[0,T_3)$
given that $T_3 = t$ is equal to the unconditioned distribution of $W^G$ during $[0,t)$. For this reason, it suffices for (\ref{eqcoup}) to show that, for each $t \geq 0$,  
\begin{equation}\label{eqcoupnew}
 \mathbb{P} \Big(  d \big( \phi, W^G (s) \big) \leq  d \big( \phi, W^G (t) \; \textrm{for all $s \in [0,t]$} \big)  \Big) \geq 1 - 2\degg^{-1} \, .
\end{equation}
It clearly suffices to prove (\ref{eqcoupnew}) in the case that $G = \tdegg$. By further conditioning on the number of transitions made by $W^{\tdegg}$ during $[0,t]$, it is enough to argue that, for each $m \in \N$,
the jump-chain $J$ of the process $\big\{ d \big(\phi, W^{\tdegg}(t) \big) : t \geq 0 \big\}$
satisfies 
\begin{equation}\label{eqcoupnewtwo}
\mathbb{P} \Big( J(n) \geq J(m) \, \textrm{for all $0 \leq m \leq n$} \Big) \geq 1 - 2 \degg^{-1}
\end{equation}
for each $n \in \N$. For given $n \in \N$, the process $R$ given by  $R(i) = J(n- i) - J(n)$ and stopped on hitting $-J(n)$
has the law of the nearest-neighbour random walk that makes a transition to the left with probability $1 - \degg^{-1}$ (also stopped on hitting $-J(n)$). The process $R$ thus remains positive at all positive times 
with probability at least $1 - 2 \degg^{-1}$. This establishes (\ref{eqcoupnewtwo}), and thus (\ref{eqcoupnew}), (\ref{eqcoup}) and the statement of the lemma. \hfill $\Box$
\noindent{\bf Proof of Lemma \ref{lemapptwo}.}
This follows directly from the proof of Lemma \ref{lemrapidadvance} in light of Lemma \ref{lemappone}. 
Indeed, the proof of Lemma \ref{lemrapidadvance} identifies an event of probability $\frac{\degg - 1}{\degg + 1}$ whose occurrence ensures that the events in the statement of the lemma occur simultaneously. \hfill $\Box$ \\
\noindent{\bf Proof of Lemma \ref{lemappthree}.} This coincides with that of Lemma \ref{lemsum}. \qed \\
\noindent{\bf Proof of Theorem~\ref{thmtwo}.}
To establish that Proposition~\ref{propone} is valid for the values of $\degg$ and $T$ asserted in Theorem~\ref{thmtwo}, we reprise the proof of the proposition. In addition to the hypothesis that $d \geq 55$  and $T \geq 101 \degg^{-1}$
(which permits us to apply the claims stated at the beginning of Appendix $A$ in place of their counterparts in the earlier proof), we must check the following analogue of~(\ref{eqcdtcb}):
$$
  \ctbetastar \cbetastar T - 2 \big( 1 - \ctbetastar \big) > 0 \, .
$$
For given $\degg$, the left-hand side is non-decreasing in $T$; for $\degg \geq 55$, it is positive when $T = 101 \degg^{-1}$. 
This completes the proof of Theorem~\ref{thmtwo}. \qed

\section{A quantitative version of Angel's result}
Here we prove Theorem~\ref{thmthree}. It follows immediately from the next result.
\begin{theorem}\label{thmappb}
Let $\de \geq 2$. Let $\tree_{[\de]}$ denote the regular rooted tree each of whose vertices except $\phi$ has $\de$ neighbours. 
If $\de \geq 765$ then $\big[ \dei + 3 \deit, \infty) \subseteq \tint^{\tree_{[\de]}}$.
For any $\epsilon > 0$, there exists $d_0' \in \N$ such that if $\de \geq d_0'$
then  $\big[ \dei + (\tfrac{13}{6} + \epsilon) \deit, \infty) \subseteq \tint^{\tree_{[\de]}}$.
\end{theorem}
The ingredients needed to prove Theorem \ref{thmappb} are Theorem~\ref{thmtwo} and the results of \cite{angel}, and the content of Appendix $B$ is a quantitative review of some of Angel's arguments. 
Note that Angel uses the symbol $\degg$ to refer to the vertex-degree of a regular tree, which corresponds to $\degg + 1$ in our usage. In the appendix, we write $\de = \degg + 1$, so that $\de$ should be identified with $\degg$ for the reader who reviews the working alongside \cite{angel}. \\ 
\noindent{\bf Proof of Theorem~\ref{thmappb}.} 
Lemma~\ref{lemcomp} implies the first statement.  The second statement of \cite[Theorem 3]{angel} and 
Lemma~\ref{lemcomp}(3) imply the second.  \qed
\begin{lemma}\label{lemcomp}
 $\empty$ 
 \begin{enumerate}
 \item $\de \geq 56$ and $\dei + 3 \deit \leq T \leq 2 \dei$ implies  $T \in \tint^{\tree_{[\de]}}$. 
 \item $\de \geq 765$ and $2 \dei \leq T \leq 101 \dei$ implies  $T \in \tint^{\tree_{[\de]}}$.
 \item $\de \geq 56$ and $T \geq 101 \dei$ implies  $T \in \tint^{\tree_{[\de]}}$.
 \end{enumerate}
\end{lemma}
\noindent{\bf Remark.}
It is noted in the proof of  \cite[Theorem 3]{angel} that a sufficient condition for  $T \in \tint^{\tree_{[\de]}}$ is 
\begin{equation}\label{expr}
(\de - 1) e^{-T}  \int_0^T \Big( e^{-a} + e^{-(T-a)} - e^{-T} \big)^{\de - 2} \dd a > 1 \, .
\end{equation}
In the derivation of this fact on page $15$, a typographical change is needed: in the second paragraph on that page, the probability that the edge $(u,v)$ rings exactly once is $Te^{-T}$, rather than the stated $e^{-T}$. However, the extra factor of $T$ is cancelled by the normalization $T^{-1}$ for Lebesgue measure on $[0,T]$, so that the derivation and form of (\ref{expr}) is not otherwise changed. \\
\noindent{\bf Proof of Lemma \ref{lemcomp}: (1).}  
As  \cite{angel} notes,
$$
 e^{-a} + e^{-(T-a)} - e^{-T}  \geq 1 - a(T - a) - T^3/6 \, .
$$
The integrand in (\ref{expr}) is at least $1 - (\de - 2)a(T-a) - (d-2)T^3/6$; thus, the integral is at least 
$T - \tfrac{1}{2}(\de - 2)T^3 - \tfrac{1}{6}(\de - 2)T^4$. Using the power series of $e^{-T}$ to cubic order as a lower bound, and the two bounds $T \geq \dei$ and $T \leq 2 \dei$ to simplify, we find that (\ref{expr}) is at least
\begin{equation}\label{exprone}
 \big(\de - 1 \big) \Big( 1 - T + \deit/2 - 4\deith/3 \Big)\Big( T - \de T^3/6 - 8\deith/3 \Big) \, .
\end{equation}
By expanding the products, omitting all positive terms whose power in $\dei$ is at least three, and using the bound $\degg \geq \gamma$ (i.e., $T \leq 2 \dei$) to estimate all negative terms of such a power, we find that
(\ref{exprone}) is at least 
$$
 1 + \big( \gamma - 13/6 \big) \dei + \big( 31\gamma/6 + 395/36 \big) \deit \, , 
$$
provided that $\de \geq \gamma$.
This expression exceeds one if $\de > f(\gamma)$, where
$$
f(\gamma)  = \frac{31\gamma/6 + 395/36}{ \gamma - 13/6} \, .
$$
Noting that $f$ is decreasing at values exceeding $13/6$ and that $f(3) \leq 32$ implies that $\de \geq 32$ and $3 \leq \gamma \leq \de$ (i.e., $\dei + 3 \deit \leq T \leq 2 \deit$) implies  $T \in \tint^{\tree_{[\de]}}$. This proves the first part of the lemma.\\
\noindent{\bf Proof: (2).}
Angel also notes that the left-hand side of (\ref{expr}) is bounded above by $2\big( e^{-T} - e^{-T\de/2} \big)$.
If $T \leq 101 \dei$ then $e^{-T} \geq 1 - 101 \dei$. Thus a sufficient condition for (\ref{expr})
is  $1 - 101 \dei > 1/2 + e^{-1}$, or $\de > \tfrac{202e}{e-2}$;  $\de \geq 765$ is sufficient for this. Hence, the second part of the lemma.\\ 
\noindent{\bf Proof: (3).} This is a restatement of Theorem~\ref{thmtwo}. \qed

\begin{lemma}
  $\de \geq 2544$ and $2 \dei \leq T \leq 0.14$ implies  $T \in \tint^{\tree_{[\de]}}$.
\end{lemma}
\noindent{\bf Proof.} Angel also notes that the left-hand side of (\ref{expr}) is bounded above by $2\big( e^{-T} - e^{-T\de/2} \big)$.
If $T \geq 2 \dei$ then $e^{-\de T/2} \leq e^{-1} < 0.368$,
while $e^{-T} > 0.868$ if $T \leq 0.14$.  \qed

\bibliographystyle{plain}

\bibliography{treecyclesbib}

\end{document}